\documentclass[11pt]{article}
\usepackage{amsfonts,amssymb,amsmath}
\usepackage{amsfonts,amssymb,amsmath,bbm,hyperref}
\hsize \textwidth
\advance \hsize by -\marginparwidth
\evensidemargin \oddsidemargin
\usepackage{amssymb}
\advance\hoffset by 5mm
\newcommand{\pdr}[2]{\frac{\partial{#1}}{\partial{#2}}}
\newcommand{\Rm}{{\mathbb R}}

\def\R{\mathbb{R}}

\def\epsilon{\varepsilon}
\def\tilde{\widetilde}
 
\newcommand{\commentout}[1]{}

\def\LL{{\mathcal{L}}}
\def\NN{{\mathcal{N}}}
\def\M{{\cal{M}}}

\def\un{{\mathbf{1}}}

\newcommand{\farc}{\frac}
\newcommand{\no}{\nonumber}

\newcommand{\norm}[1]{\lVert #1 \rVert}

\newcommand{\br}{\begin{eqnarray}}
\newcommand{\er}{\end{eqnarray}}
\newcommand{\be}{\begin{equation}}
\newcommand{\ee}{\end{equation}}
\newcommand{\baa}{\begin{array}}
\newcommand{\eaa}{\end{array}}
\newcommand{\ba}{\begin{eqnarray}}
\newcommand{\ea}{\end{eqnarray}}
\def\di{\displaystyle}
\def\ln{{\log}}

\numberwithin{equation}{section}

\newtheorem{theorem}{\bf Theorem}[section]

\newtheorem{thm}[theorem]{Theorem}

\newtheorem{cor}[theorem]{Corollary} 
 
\begin{document}

\title{Refined long time asymptotics for Fisher-KPP fronts} 

\author{ James Nolen\thanks{Department of Mathematics, 
Duke University, Durham, NC 27708, USA; nolen@math.duke.edu} \and 
Jean-Michel Roquejoffre\thanks{Institut de Math\'ematiques (UMR CNRS 5219),
Universit\'e Paul Sabatier, 118 route de Narbonne, 31062 Toulouse
cedex, France; jean-michel.roquejoffre@math.univ-toulouse.fr} \and Lenya
Ryzhik\footnote{Department of Mathematics, Stanford University,
Stanford CA, 94305, USA; ryzhik@math.stanford.edu } }

\date{February 23, 2018}
\maketitle


\begin{abstract}
We study the one-dimensional Fisher-KPP equation, with an initial condition $u_0(x)$ that coincides with the step function
except  on a compact set. A well-known result of 
M. Bramson in~\cite{Bramson1,Bramson2} states that, as $t\to+\infty$,
the solution   converges to a traveling wave located at the position~$X(t)=2t-(3/2)\log t+x_0+o(1)$,
with the shift $x_0$ that depends on~$u_0$.
U.~Ebert and W. Van Saarloos have formally derived in~\cite{EVS,VS} a correction to the Bramson shift, 
arguing that~$X(t)=2t-(3/2)\log t+x_0-3\sqrt{\pi}/\sqrt{t}+O(1/t)$. Here, we prove that this result does hold, with an error
term of the size $O(1/t^{1-\gamma})$, for any $\gamma>0$. The interesting aspect of this asymptotics 
is that the coefficient in front of the
$1/\sqrt{t}$-term does not depend on $u_0$. 
 
\end{abstract}

\section{Introduction}

The goal of this paper is to provide a sharp  
large time asymptotics of the solutions   the Fisher-KPP equation 
\begin{equation}
\label{e2.1}
u_t-u_{xx}=u-u^2,~~t>0,~x\in\Rm.
\end{equation}
The initial condition $u_{in}(x)=u(0,x)$ is a compactly supported perturbation of the step function:
there exists $L>0$ so that $u_{in}(x)\equiv 1$ for $x<-L$ and $u_{in}(x)\equiv 0$ for $x\ge L$. In addition,
we assume that $0\le u_{in}(x)\le 1$ for all $x\in\Rm$, so that $0<u(t,x)<1$ for all $t>0$ and $x\in\Rm$.
The assumptions on the initial condition, especially as $x\to-\infty$ can be significantly weakened,
without any change in the result. The more stringent conditions are adopted purely for convenience, but 
we stress that  the decay of $u_0(x)$ as $x\to+\infty$ does have to be faster than $\exp(-x)$ for the results
to hold. For a detailed study of this issue we refer to~\cite{BBHR} where a related linear problem with similar
properties has been studied. 

This issue has a long history. The first contribution is that of Fisher \cite{Fisher}, who
identified the spreading velocity~$c_*=2$ of the solutions
via numerical computations and other arguments. 
In the same year, the pioneering KPP paper \cite{KPP} proved 
that the solution of~\eqref{e2.1}, starting from 
a step function, converges to a traveling wave profile in the following sense: there is a   function
\[
\sigma_\infty(t)=2t+o(t),~~\hbox{ as $t\to+\infty$},
\]
such that
\begin{equation}\label{sep1820}
\lim_{t\to+\infty}u(t,x + \sigma_\infty(t))=\phi(x).
\end{equation}
Here, $\phi(x)$ is the profile of a  traveling wave 
that connects the stable equilibrium~$u\equiv 1$ to the unstable equilibrium~$u\equiv 0$ and
moves with the minimal speed $c_*=2$:
\begin{equation}
\label{e2.4}
\begin{array}{rll}
&-\phi''-2\phi'=\phi-\phi^2,  \\
&\phi(-\infty)=1, \quad \phi(+\infty)=0.
\end{array}
\end{equation}
Each solution $\phi(\xi)$ of 
\eqref{e2.4} is a shift of a fixed profile $\phi_*(\xi)$: 
$
\phi(\xi)=\phi_*(\xi+s),
$
with some fixed $s\in\Rm$. The function $\phi_*(\xi)$
has the asymptotics
 \begin{equation}
\label{e2.20}
\phi_*(\xi)=(\xi+k)e^{-\xi}+O(e^{-(1+\omega_0 )\xi}),
\end{equation}
with two universal constants $\omega_0 >0$, $k\in\R$.
The question whether the function $\sigma_\infty(t)$ tends to a constant, or is a 
nontrivial sublinear function of time, was solved by   Bramson
\cite{Bramson1}, \cite{Bramson2}.  
\begin{thm}
\label{t2.10}\cite{Bramson1,Bramson2}
There is  a constant $x_\infty$, depending on the initial condition $u_0(x)$, such that
\begin{equation}\label{apr604}
u(t,x)=\phi_*(x-2t+\frac32\log t-x_\infty)+o(1),\hbox{ as {$t\to+\infty$}},
\end{equation}
in the sense of uniform convergence on $\Rm$.
\end{thm}
Both papers by Bramson use probabilistic tools, and elaborate explicit computations. 
The reason why the probabilistic arguments are natural here is that \eqref{e2.1} is  
related to the branching Brownian motion~\cite{McK}. This connection brought a lot of 
recent activity on the Fisher-KPP equation in the probability and
physics communities -- see, for instance, \cite{BD2,BD1}.
The results of~\cite{Bramson1,Bramson2}
were also proved by Lau~\cite{Lau}, 
using the decrease of  the number of intersection points 
between any two solutions of the parabolic Cauchy problem~(\ref{e2.1}).  

A short and simple proof of  Theorem \ref{t2.10}, solely relying on the PDE arguments,
was given recently in \cite{HNRR1,NRR1}: first,
the estimate 
\[
\sigma_\infty(t)=2t-\di\frac32\log t+O(1)
\]
was proved in~\cite{HNRR1}, and then the 
full estimate 
\begin{equation}\label{apr602}
\sigma_\infty=2t-\di\frac32\log t+x_\infty,
\end{equation}
with $x_\infty$ depending on the initial datum, was proved in~\cite{NRR1}.
The ideas of \cite{HNRR1} were developed in a more complex paper \cite{HNRR2}
to compute a logarithmic shift in a version of \eqref{e2.1} with spatially periodic coefficients, 
a situation that had not been treated previously by the probabilistic methods. 

The $\log t$ correction in (\ref{apr602}) is unusual: for  reaction-diffusion equations of the type
\[
u_t-u_{xx}=f(u),\ \ t>0,\ x\in\R
\]
one sees, most of the time,  exponential in time convergence  to a
constant shift of a traveling wave, see for instance the classical Fife-McLeod paper \cite{FML}. 
This raises the question of the convergence rate in \eqref{apr604}.
That is, the issue is to estimate the error between  
\begin{equation}\label{oct730}
\sigma(t)=\sup\{x:~u(t,x)=1/2\}, \hbox{  and}\ \ \ \bar\sigma_\infty(t):=2t-\di\frac32{\mathrm{log}}t+x_\infty.
\end{equation}
A very interesting paper of Ebert and Van Saarloos~\cite{EVS}, completed in~\cite{VS}, performs a formal
analysis of the convergence  and states that 
\begin{equation}\label{sep1806} 
\sigma(t)=\bar\sigma_\infty(t)-\farc{3\sqrt\pi}{\sqrt{t}}+o(\frac1{\sqrt t}).
\end{equation}
A striking feature is that the predicted constant  $3\sqrt{\pi}$ in (\ref{sep1806}) 
does not depend on the 
initial condition, unlike the zero order term $x_\infty$.

Here, we prove a rigorous version of (\ref{sep1806}). We do this by constructing
an approximate solution of~(\ref{e2.1}), which is approached by the solutions of \eqref{e2.1} at a rate almost equal to~$O(t^{-1})$. Examination 
of the shift of the approximate solution provides the asymptotics of $\sigma(t)$.

\subsection*{Main results}

One of the main ingredients in this paper is the construction of an approximate solution which solves the equation up to  a
sufficiently small correction. Here is the precise result.
\begin{thm}
\label{t2.1} For all $\gamma\in(0,1/10)$,
there is a one-parameter family $(u_{app}(t,x+\lambda))_{\lambda\in\R}$ of   the form
\[
u_{app}(t,x)=\phi_*(x-\tilde\sigma(t))+u_0(t,x-\tilde\sigma(t))+\frac{u_1(t,x-\tilde\sigma(t))}{\sqrt t},
\]
with
\begin{equation}\label{apr702}
\tilde\sigma(t)=2t-\di\frac32{\log}t-\di\frac{3\sqrt\pi}{\sqrt t}+O(\frac1{t^{1-\gamma}}).
\end{equation}
The functions $u_0(t,x)$ and $u_1(t,x)$ are bounded and continuous, and supported in $\{x>t^\gamma\}$.
In addition, $u_0$ is of the class $C^1$, and $u_1$ is $C^1$ everywhere except at 
$x=t^\gamma$, where it has a jump of the~$x$-derivative.  The functions $u_{app}(t,x)$ are approximate solutions to (\ref{e2.1}) 
in the sense that
\begin{eqnarray}
\label{e1.30}
&& \Big\vert\big(\partial_tu_{app}-\partial_{xx}^2u_{app}-u_{app}+u_{app}^2\big)(t,x+\tilde\sigma(t))\Big\vert
\le  C_\gamma t^{-1+2\gamma} \left(e^{-x} \un_{0<x<t^\gamma}+\un_{x< 0} \right)\\
&&~~~~~~~~~~~~~~~~~~~~~~~~+ C_\gamma t^{-3/2} e^{-x-x^2/((4+\gamma)t)}\un_{x>t^\gamma}
+ C_\gamma t^{-1+2\gamma}\delta(x-t^\gamma).\nonumber
\end{eqnarray}
\end{thm}
The estimate in the right side  includes the spatial 
behavior of the error -- this is needed in the region
where the solution is small. The different error sizes in the regions $x<t^\gamma$ and $x>t^\gamma$ in
(\ref{e1.30}) come about because we need less precision in approximating 
the solution to the left of $x=t^\gamma$, where $u$ is either $O(1)$ or not too 
small, than to the right of $x=t^\gamma$, where $u$ is ``very small".
The delta function in the last term in the right side is not an issue, and can be, 
in principle, eliminated by a modification of the
approximate solution. 
With this result in hand, the next task is to prove that the solutions of \eqref{e2.1} converge to 
a shift of $u_{app}$ at a certain rate.  Our second main result is:
\begin{thm}
\label{t2.2}
For all $\gamma>0$, there is $C_\gamma>0$ such that, 
for all $t\geq 0$ and all $x\in\Rm$, we have, with~$\tilde\sigma(t)$
as in (\ref{apr702}), and some $x_\infty\in\Rm$, depending on the initial
condition $u_{in}$:
\begin{equation}
\label{e1.150}
\vert u(t,x+\tilde\sigma(t))-
u_{app}(t,x+\tilde\sigma(t)+x_\infty)\vert
\leq\frac{C_\gamma (1+\vert x\vert)e^{-x}}{t^{1-\gamma}}.
\end{equation}
\end{thm}
The corollary of this result is the following
\begin{cor}
\label{c2.2}
If we fix $s\in(0,1)$ and define the front position as
$
\sigma_s(t)=\max\{x:~u(t,x)=s\},
$
then $\sigma_s(t)$ has an asymptotics of the form 
\[
\sigma_s(t)=2t-\di\frac{3}{2}\log t+x_{\infty} + \phi_*^{-1}(s) - \di\frac{3\sqrt\pi}{\sqrt{t}}+O(\farc{1}{t^{1-\gamma}}).
\]
\end{cor}
This confirms the Ebert-Van Saarloos prediction.

\subsubsection*{Related works}
 
The $3\sqrt\pi$  
prediction  has already been verified by C.~Henderson in~\cite{Henderson},
for a linearized moving boundary problem:
\begin{eqnarray}
\label{e5.200}
&&U_t-U_{xx}=U,~~ t>0,x>\sigma(t),\\
&&U(t,\sigma(t))=0,\nonumber
\end{eqnarray}
and a compactly supported initial condition. 
The Dirichlet boundary condition serves the same purpose
as the term $(-u^2)$ in the KPP equation -- when the moving boundary is chosen ``correctly", the solution
of (\ref{e5.200}) does not grow or decay  in time. Both solutions of~(\ref{e2.1}) and~(\ref{e5.200})
are governed by the ``far ahead" tails where they are small -- these are so called pulled fronts. The
difference between (\ref{e5.200}) and the full
KPP problem on the whole line is that (\ref{e2.1}) has an ``inner" layer where the solution transitions
from $O(1)$ to very small values. The moving boundary in \cite{Henderson}
is taken of the form
\[
\sigma(t)=2t-\di\frac32\log t-\di\frac{c}{\sqrt{t}},
\]
for $t\ge 1.$
Then, if $c=3\sqrt\pi$, there is $\alpha_0>0$ such that 
\begin{equation}\label{apr704}
\biggl\vert \int_{\sigma(t)}^{+\infty}U(t,x)dx-\alpha_0 \biggl\vert
\leq\frac{C\log t}t.
\end{equation}
On the other hand, if $c\neq3\sqrt\pi$, the convergence rate 
in (\ref{apr704}) is of the order $1/ \sqrt{t}$. 
We refer to a recent preprint  \cite{BBHR} for a very 
detailed study of the same problem, 
according to the behavior of the initial condition at infinity.  

As we have mentioned, an interesting feature of the problem is that the $t^{-1/2}$ 
correction to the Bramson shift is universal, in the sense that it is  
independent of the initial datum. In addition, the analysis can be easily adapted 
to show that an identical result holds for more general equations of the
form
\[
u_t=u_{xx}+f(u),
\]
with a KPP type nonlinearity: $f\in C^1[0,1]$, $f(0)=f(1)=0$, and $f(u)\le f'(0)u$ fot all $u\in(0,1)$.
In that case, the ``$3\sqrt{\pi}/\sqrt{t}$" term in the
shift  depends on the nonlinearity $f(u)$ only through
$f'(0)$, and the shape of the
solution approaches the traveling wave profile at a rate almost $O(t^{-1})$. The 
preprint~\cite{BB} explains why this last  feature holds: if 
the $t^{-1/2}$ correction were to depend of the value of the solution, this would entail 
wild oscillations to the front, that are not confirmed by the numerics. This result 
was a strong incentive for us to verify the actual value of the coefficient in front of $1/\sqrt{t}$.

{\bf Organization of the paper.}
In Section~\ref{sec:2}, we explain, in an informal way, why the results  are likely 
to hold.  We then prove Theorem \ref{t2.1}
in Section~\ref{sec:approx}, where  we construct the approximate solution.   In 
Section~\ref{sec:true}, we use the approximate solution to prove Theorem \ref{t2.2} 
and its corollaries.

{\bf Acknowledgment.} JN was  supported by NSF grant DMS-1351653, 
and LR by NSF grant DMS-1311903. JMR  was supported by  the European Union's Seventh 
Framework Programme (FP/2007-2013) / ERC Grant
Agreement n. 321186 - ReaDi - ``Reaction-Diffusion Equations, Propagation and 
Modelling'',
as well as the ANR project  NONLOCAL ANR-14-CE25-0013.  
LR and JMR thank  the Labex CIMI for  a 
PDE-probability  quarter in Toulouse, in Winter 2014, out of which the idea of this
paper grew 
and which provided a stimulating scientific environment for this project. The 
authors thank H. Berestycki for raising the issue of a mathematically rigorous proof 
of \cite{EVS}. Finally, we thank E. Brunet and J. Berestycki 
for educating us about the problem, pointing out a computational
mistake in an earlier version of this paper, and for making earlier versions of 
their preprint \cite{BB} available to us.

\section{Strategy of the proofs}\label{sec:2}

 Consider the Cauchy problem (\ref{e2.1}) 
starting at $t=1$ for convenience of the notation:
\begin{eqnarray}\label{sep2704}
&&u_t-u_{xx}=u-u^2,~~x\in\Rm,~~t>1,
\\
&&u(1,x)=u_{in}(x)=1-H(x)+v_0(x),\  \ \hbox{$v_0$ compactly supported},\nonumber
\end{eqnarray} 
 and
proceed with the standard sequence of changes of variables
\begin{equation}\label{apr706}
x \mapsto x - 2t + ({3}/{2}) \log t,\quad u(t,x)=e^{-x}v(t,x)
\end{equation}
so that $v$ solves
\begin{equation}
\label{e2.6}
v_t-v_{xx}-\frac3{2t}(v-v_x)+e^{-x}v^2=0, \quad \quad x \in \Rm, \quad t > 1.
\end{equation}

We stress that 
the removal of the exponential factor in (\ref{apr706}) is critical for
understanding the dynamics of $u(t,x)$ as ``basically diffusive".  

For any $x_\infty \in \Rm$, the function 
\[
V(x) = e^{x} \phi(x - x_\infty)
\]
satisfies 
\begin{equation}\label{apr708}
V_t - V_{xx} + e^{-x} V^2 = 0.
\end{equation}
Note that (\ref{e2.6}) is a perturbation of (\ref{apr708})
for $t\gg 1$, and both of them are close to the diffusion equation for~$x\gg 1$.  
Hence, ``everything" relevant to the solutions 
of ~\eqref{e2.6} should happen at the diffusive spatial scale~$x\sim\sqrt t$. It is convenient to pass to the self-similar variables
\begin{equation}
\label{e2.500}
\tau=\log t,\ \ \ \eta=\frac{x}{\sqrt{t}}.
\end{equation}
This transforms \eqref{e2.6} into
\begin{equation}
\label{e2.7}
w_\tau-\di\frac\eta2w_\eta-w_{\eta\eta}-\di\frac32w
+\di\frac32e^{-\tau/2}w_\eta+e^{\tau-\eta{\mathrm{exp}}(\tau/2)}w^2=0,\  \   \   \eta\in\R,~~\tau>0.
\end{equation}
It is easy to see now why the linearized problem with the Dirichlet boundary
condition at $\eta=0$ is a good approximation to (\ref{e2.7}). Indeed,
for $\eta<0$, the last term in the left side of (\ref{e2.7}) becomes very large,
which forces $w$ to be very small in this region. 
On the other hand, for~$\eta>0$, this term is very small, so it 
should not play any role in the dynamics of $w$ for $\eta>0$. The main step in the argument of~\cite{NRR1} (see Lemma 5.1 therein) is
a convergence result of the form
\begin{equation}\label{sep2702}
w(\tau,\eta)\sim \alpha_\infty \eta e^{\tau/2-\eta^2/4},~~\eta>0.
\end{equation}
More specifically, as $\tau \to \infty$, $e^{-\tau/2} w(\tau,\eta)$ converges in $L^2(0,\infty)$ to $\alpha_\infty \eta e^{-\eta^2/4}$. 
Therefore, we have (reverting to the variables of (\ref{e2.6}))
\begin{equation}\label{nov2302}
u(t,x)=e^{-x}v(t,x)\sim \alpha_\infty xe^{-x}e^{-x^2/(4t)},
\end{equation}
at least for $x$ of the order $O(\sqrt{t})$.  
This, in view of the asymptotics (\ref{e2.20}) of the wave $\phi_*$ at infinity
determines the unique translation: 
\begin{equation}
\label{e2.200}
x_\infty= \log\alpha_\infty.
\end{equation}

This argument gives the right insight for the construction of the approximate solution. The idea is to view $1/{\sqrt t}$ as a 
small parameter, in terms of which one may expand the solution. 
It is natural to identify two zones: the 
region near the front, that is, $x\sim O(1)$ --
it corresponds to~$\eta\sim e^{-\tau/2}$, a very small region in the self-similar variables,  
and the diffusive region, where $x\sim\sqrt t$ and $\eta\sim O(1)$. 
The transition region is $x\sim t^{\gamma}$, with $\gamma>0$ 
small. We  perform a classical  asymptotic expansion of an inner solution in the region $x\sim O(1)$, 
approximating $u$ near the front, and of an outer solution, approximating 
$u$   at distances $O(\sqrt{t})$ from the front. 
Matching the inner and outer expansions is done in the intermediate 
region $x\sim t^\gamma$. 
 
Once the translate $x_\infty$ is selected, this also determines the translate of the
approximate solution to which the solution is supposed to converge, at a rate 
faster than $t^{-(1-\gamma)}$, for all
small $\gamma$. Everything reduces to proving that the difference between the true solution and the approximate solution 
will not exceed $t^{\gamma-1}$.
The argument is long and technical, and is carried out in the 
self-similar variables~\eqref{e2.500}.  However, it relies on two simple ideas. 
The first is to transform the problem on the whole line into a Dirichlet problem on 
the half line, by a classical sequence of transformations 
and the final subtraction of the value of $u$ at $t^{\gamma}$.  
The trouble is that the nonlinear term $u^2$ in the original equation~\eqref{e2.1} 
provides, as 
usual, a term which may grow like $e^{3\tau/2}$ in (\ref{e2.7}). 
The difficulty 
is overcome by noticing that its  support shrinks as $e^{-\tau/2}$. A large 
part of the proof is devoted to estimating this term in the best 
way.  For that, we first obtain weak estimates on the difference~$u-u_{app}$, which still yield an 
improvement of the nonlinear term. This improvement entails a 
better estimate on $u-u_{app}$, and so on. As we have mentioned, the technical details are 
nontrivial.
 
 \section{The approximate solution}\label{sec:approx}
 
Instead of working directly with \eqref{e2.6}, we introduce the moving frame
that incorporates a (still unknown) correction of the 
order $t^{-1/2}$, namely, instead of (\ref{apr706}), 
we make a slightly different
successive change of variables:
\[
x \mapsto x - 2t + ({3}/{2}) \log t-\frac\sigma{\sqrt t},
\quad u(t,x)=e^{-x}v(t,x).
\]
The function $v$ satisfies
\begin{equation}
\label{e4.11}
v_t-v_{xx}-(\frac3{2t}+\frac\sigma{2t^{3/2}})(v-v_x)+e^{-x}v^2=0, \quad \quad x \in \Rm, \quad t > 1.
\end{equation}
Let us denote this nonlinear operator as
\begin{equation}\label{apr712}
N\!L[v] =v_t-v_{xx}-(\frac3{2t}+\frac\sigma{2t^{3/2}})(v-v_x)+e^{-x}v^2.
\end{equation}
We will construct an approximate solution to (\ref{e4.11}), called $V_{app}(t,x)$
As we have mentioned, it is natural to 
consider an intermediate scale $x\sim O(t^\gamma)$, with some~$\gamma>0$,
and seek an approximate solution to (\ref{e4.11}) in two different forms: 
one valid for $x\leq t^\gamma$, 
the other valid for~$x\geq t^\gamma$:
\[
V_{app}(t,x)=V^-(t,x)\ \hbox{for $x<t^\gamma$,}\   \     \     \    \    \   V_{app}(t,x)=V^+(t,x)\ 
\hbox{for $x>t^\gamma$.}
\]
The functions $V^-$ and $V^+$ will be matched at $x=t^\gamma$. 
 
\subsection{The inner approximate solution $V^-$}

Note that \eqref{e4.11} contains terms that are either of order $O(1)$, or 
of the order $O(t^{-1})$ and smaller. 
So, a natural first guess is to choose $V^-(t,x)=V^-(x)$ and to  discard the $O(t^{-1})$ terms. In other words, we impose
\[
-(V^-)''+e^{-x}(V^-)^2=0.
\]
A first choice is 
\begin{equation}\label{apr716}
V_0^-(x)=e^{x}\phi_*(x).
\end{equation}
This function has the asymptotics:
\begin{equation}\label{apr720}
V_0^-(x)\sim e^x \hbox{ as
$x\to-\infty$, and $V_0(x)\sim x$ as $x\to+\infty$.}
\end{equation}
We will have to correct it slightly at $x\sim t^\gamma$ in order to ensure the 
matching with $V^+(t,x)$. Hence, we choose $V^-$ as
\begin{equation}
\label{e4.1000}
V^-(t,x)=V_0^-(x+\zeta(t))=e^{x+\zeta(t)}\phi_*(x+\zeta(t)).
\end{equation}
Here, the correction $\zeta(t)$, which will come from the matching procedure, will be of the 
order 
\begin{equation}\label{apr718}
\zeta(t)\sim O(t^{-1+3\gamma}),~~\dot\zeta(t)\sim O(t^{-2+3\gamma}).
\end{equation}

Let us now estimate $N\!L[V^-]$: 
\begin{eqnarray}\label{apr714}
&&N\!L[V^-]=\dot\zeta V_0^-(x+\zeta(t))
-(V_{0}^{-})''(x+\zeta(t))-(\frac3{2t}+\frac\sigma{2t^{3/2}})
(V_0^-(x+\zeta(t)) -(V_0^-)'(x+\zeta(t)))\nonumber\\
&&+
e^{-x}(V_0^-)^2(x+\zeta(t))=
\dot\zeta V_0^-(x+\zeta(t))
-(\frac3{2t}+\frac\sigma{2t^{3/2}})
(V_0^-(x+\zeta(t)) -(V_0^-)'(x+\zeta(t)))\nonumber\\
&&+
\Big[e^{-x}-e^{-x-\zeta(t)}\Big](V_0^-)^2(x+\zeta(t)).
\end{eqnarray}
Note that all terms in (\ref{apr714}), decay as
$e^{x}$ for $x<0$ because of (\ref{apr720}).
Taking also into account (\ref{apr718}) gives
\begin{equation}
\label{e4.24}
N\!L[V^-](t,x)=
n_1(t,x)(\un_{0<x<2t^\gamma}(x)+\un_{\R_-}(x)e^x), \quad \quad x \leq 2t^\gamma,
\end{equation}
with
\begin{equation}\label{oct1102}
|n_1(t,x)|\le Ct^{-1+3\gamma}.
\end{equation}

\subsection{The outer approximate solution $V^+$}

In the outer region $x>t^\gamma$, we pass to the self-similar variables
\begin{equation}
\label{e4.34}
\tau=\ln \,t,\  \   \   \eta=\frac{x+x_0}{\sqrt{t}},
\end{equation}
the shift $x_0$ kept free for the moment. Our starting point is,
again, (\ref{e4.11}), in the self-similar variables. 
The equation for $V^+$ is 
\begin{equation}
\label{e4.21}
v_\tau-v_{\eta\eta}-\di\frac\eta2v_\eta+(\di\frac32+\farc{\sigma}{2} e^{-\tau/2})(e^{-\tau/2}v_\eta-v)+e^{\tau-\eta e^{\tau/2}+x_0}v^2=0.
\end{equation}
We will set
\begin{equation}
\label{e3.1001}
Lv=-v_{\eta\eta}-\di\frac{\eta}{2}v_\eta-v.
\end{equation}
As in the construction of $V_{app}^-$, we are not going to solve (\ref{e4.21}) 
exactly,
but find an approximate solution. 
Strictly speaking, we only 
need $V^+$ defined for $x>t^\gamma$,
that is, for~$\eta> e^{-(1/2 - \gamma) \tau}$ but we will define it for $\eta\ge 0$. We impose the 
boundary condition
\begin{equation}\label{oct502}
V^+(\tau,0)=0,
\end{equation}
which is consistent with the presence of the absorption term $e^{\tau-\eta e^{\tau/2}}v^2$ in the left
side of (\ref{e4.21}),
which is huge as soon as $\eta$ is just a little negative.
As $V^{-}(t,x)$ is of the order $O(t^\gamma)$ at $x=t^\gamma$, to have
a hope of a good matching we need 
\[
V^+(\tau,e^{-(1/2-\gamma)\tau})\sim e^{\gamma\tau}.
\]
On the other hand, the boundary condition (\ref{oct502}) means that
\[
V^+(\tau,e^{-(1/2-\gamma)\tau})\sim \pdr{V^+(\tau,0)}{\eta}e^{-(1/2-\gamma)\tau},
\]
thus we need
\[
\pdr{V^+(\tau,0)}{\eta}\sim e^{\tau/2}.
\]
Hence, it is natural to look for $V^+$ in the form 
\[
V^+(\tau,\eta)=e^{\tau/2}V_0^+(\eta)+V_1^+(\eta).
\]
Inserting this ansatz into (\ref{e4.21}) and collecting the leading order terms gives
\begin{equation}\label{oct506}
LV_0^+=0,
\end{equation}
and 
\begin{equation}\label{oct508}
(L-\farc{1}{2})V_1^++\di\frac32(V_0^+)_\eta-\frac\sigma2V_0^+=0,
\end{equation}
with the boundary conditions 
\begin{equation}\label{oct510}
V_i^+(0)=V^+_i(+\infty)=0,~~i=0,1.
\end{equation}
Setting 
\[
e_0(\eta)=\eta e^{-\eta^2/4}\hbox{ for $\eta>0$}, 
\]
we have 
\begin{equation}\label{oct512}
V_0^+(\eta)=q^+_0e_0(\eta),
\end{equation}
the constant $q^+_0$ being for the moment free. Once $V_0^+$ is fixed, there is a unique solution $V_1^+$  to~(\ref{oct508}), 
with $e^{\eta^2/(4 + \gamma)} V_1 \in L^2(\Rm_+)$, because the
spectrum of  $L$ is $\{0,1,2,\dots\}$.

We will need the derivative $(V_1^+)_\eta(0)$ for the matching procedure. 
The (formal) adjoint of $L$ satisfies
\begin{equation}
\label{e4.1003}
L^*(1-\frac{\eta^2}2)=0.
\end{equation}
Multiplying \eqref{oct508} by $1-\eta^2/2$ and integrating by parts gives
\begin{equation}
\label{e4.1002}
(V_1^+)'(0)=\int_0^{+\infty}(1-\frac{\eta^2}2)(\frac\sigma2V_0^+-\frac32(V_0^+)')d\eta=-[\sigma+3\sqrt\pi]q_0^+.
\end{equation}

\subsubsection*{Estimating the error}

Let us denote by $\NN\!\LL[v]$ the nonlinear operator in the left side of \eqref{e4.21}. Then we have
\begin{equation}
\label{e4.25}
|\NN\!\LL[V^+]|\le Ce^{-\tau/2}
\un_{\R_+}(\eta)e^{-\eta^2/(4+\gamma)}.
\end{equation}
In the original variables, the function $V^+$ has the form
\begin{equation}\label{oct702}
V^+(t,x)=q_0^+(x+x_0)e^{-(x+x_0)^2/(4t)}+V_1^+\Big(\frac{x+x_0}{\sqrt{t}}\Big),
\end{equation}
and (\ref{e4.25}) implies that
\begin{equation}\label{oct718}
|N\!L[V^+](t,x)|\le Ct^{-3/2}\un_{\{x + x_0 > 0\}}e^{-(x+x_0)^2/((4+\gamma)t)}, \quad \quad \text{for}\;\;x \geq -x_0.
\end{equation}
Here, $N\!L[V^+]$ is as in (\ref{apr712}).

\subsection{Matching the inner and outer approximate solutions} 

Our next task is to choose the parameters so that the inner and outer approximate solutions match at $x=t^\gamma$. Ideally, we would like to match 
both $V^-$ and $V^+$ and their derivatives at this point. 
However, $V^-$ and $V^+$ are of the size $O(t^\gamma)$ in 
this region -- they are ``large", while their derivatives are $O(1)$.
Thus, the key is to match $V^-$ and $V^+$ and
the matching of the derivatives is less of an issue. 

Recall that we have
\begin{equation}
\label{e4.27}
V^-(t,t^\gamma)=t^\gamma+k+\zeta(t)+O(e^{-\omega_0t^\gamma})
\end{equation}
while for $V^+(t,t^\gamma)$, using expression (\ref{oct702}) we get
\begin{eqnarray}
\label{e4.27bis}
&&V^+(t,t^\gamma)=t^{1/2}V_0^+\Big(\farc{t^\gamma+x_0}{\sqrt{t}}\Big)+
V_1^+\Big(\farc{t^\gamma+x_0}{\sqrt{t}}\Big)\\
&&~~~~~~~~~~~=
q_0^+\biggl((t^\gamma+x_0)(1+O(t^{2\gamma-1}))
-(\sigma+3\sqrt\pi) t^{-1/2}(t^\gamma+x_0)\biggl)+O(\frac1{t^{1-2\gamma}}). \nonumber
\end{eqnarray}
Equating the terms of the order $O(t^\gamma)$ and $O(1)$ gives
\begin{equation}\label{oct610}
q_0^+=1,~~x_0=k,
\end{equation}
while those of the order $O(t^{-1/2+\gamma})$ and $O(t^{-1/2})$ give 
\begin{equation}\label{oct612}
\sigma=-3\sqrt{\pi}.
\end{equation}
Finally, we choose $\zeta(t)$ to eliminate the terms of the order higher than $O(t^{-1/2})$, which means that
\begin{equation}\label{oct614}
\zeta(t)=O(\farc{1}{t^{1-3\gamma}}).
\end{equation}
This implies, by inspection, that
\[
\dot\zeta(t)=O(\farc{1}{t^{2-3\gamma}}).
\]
Therefore, both conditions in (\ref{apr718}) are satisfied.

Choosing the parameters in this way, we have matched 
the values of $V^+$ and $V^-$ at $x = t^\gamma$: 
\[
V^+(t,t^\gamma) = V^-(t,t^\gamma),
\]
but we have no freedom left in terms of the parameters to match their derivatives at this point. This is a relatively
minor inconvenience as  $N\!L[V_{app}]$ would then have a Dirac mass, of the size proportional
to the jump in the derivatives. Taking into account (\ref{oct512}) and (\ref{e4.1002}), as well as
(\ref{oct610})-(\ref{oct614}), we see 
that these derivatives are given by:
\begin{equation}\label{oct704}
V_x^+(t,t^\gamma)=e^{-(t^\gamma+k)^2/(4t)}-
\frac{(t^\gamma+k)^2}{2t}e^{-(t^\gamma+k)^2/(4t)}+
\farc{1}{\sqrt{t}}(V_1^+)'\Big(\frac{t^\gamma+k}{\sqrt{t}}\Big)
=1+O(\frac{1}{t^{1-2\gamma}}),
\end{equation}
and,  
\begin{equation}\label{oct706}
V_x^-(t,t^\gamma)=(V_0^-)'(t^\gamma + \zeta(t))=1+
O(e^{-\omega_0 t^\gamma}).
\end{equation}
We conclude that with our choice of $V^+$ and $V^-$ the jump in the derivatives is
very small:
\begin{eqnarray}\label{oct708}
V_x^+(t,t^\gamma)-V_x^-(t,t^\gamma)\sim O(\farc{1}{t^{1-2\gamma}}).
\end{eqnarray}
We could have avoided this jump by modifying slightly the approximate
solution, at the expense of even longer formulas.

\vspace{0.2in}

{\bf Summary:} The full approximate solution $V_{app}(t,x)$
for (\ref{e4.11}) is defined by
\begin{equation}\label{apr724}
V_{app}(t,x) = V^-(t,x) \un_{x < t^\gamma} + V^+(t,x) \un_{x \geq t^\gamma}.
\end{equation}
The inner and outer pieces have the form:
\begin{equation}\label{apr726}
V^-(t,x) = e^{x+\zeta(t)} \phi_*(x+\zeta(t)),\quad\quad \zeta(t)=O(t^{3\gamma-1}),\quad \dot\zeta(t)=O(t^{3\gamma-2}),
\end{equation}
and
\begin{equation}\label{apr728}
V^+(t,x)=(x+k)e^{-(x+k)^2/(4t)}+V_1^+\Big(\frac{x+k}{\sqrt{t}}\Big),
\end{equation}
The function $V^+$ does not depend on the choice of $\gamma$, while 
$V^-$ depends on $\gamma$, through the shift~$\zeta(t)$.

Inserting the ansatz (\ref{apr724}) into 
\eqref{e4.11} yields, in view of~\eqref{e4.24}-\eqref{oct1102} 
and \eqref{oct718}, and taking into account that we use $V^-$ for $x<t^\gamma$
and $V^+$ for $x>t^\gamma$:
\begin{equation}
\label{e4.30}
|N\!L[V_{app}](t,x)| \le Ct^{-1+3\gamma}
(\un_{0<x<t^\gamma}+e^x \un_{x < 0} ))+ Ct^{-3/2} e^{-x^2/((4+\gamma)t)}\un_{x>t^\gamma}
+ Ct^{-1+2\gamma}\delta(x-t^\gamma).
\end{equation}
The first two terms come from $N\!L[V^-]$ and $N\!L[V^+]$,
respectively, while the singular term~$\delta(x - t^\gamma)$ comes from the jump 
(\ref{oct708}) in the derivative at the matching point $x=t^\gamma$. This estimate
is the main result of this section. 

\smallskip
\noindent{\bf Remark.} {\rm It is now clear why the $t^{-1/2}$ term in the expansion of the front location
does not depend on the initial datum, as it
is determined by a matching procedure that is itself independent of $u_0$. It is another manifestation of the role played by the diffusive 
zone $\{x\sim\sqrt t\}$, which actually drives the dynamics of the solution. Let us recall that the shift $x_\infty$ is also determined by the 
diffusive zone.}

\section{The approximate solution is an approximation to the true solution}\label{sec:true}

From \cite{NRR1} (and from \cite{Bramson1, Bramson2}), we know that there is an asymptotic shift
$x_\infty$ such that, as~$t\to+\infty$, we have $u(t,x)\to\phi_*(x - x_\infty)$ uniformly on $\R$. 
Without loss of generality, we will assume that the initial condition is such that
\[
x_\infty=0.
\]
As in Section \ref{sec:approx}, we will work in the frame moving as 
$2t-(3/2)\log t - 3 \sqrt{\pi/t}$.  If $u(t,x)$ is the solution of the Fisher-KPP equation in this moving 
frame, then the function
\[
v(t,x) = e^{x} u(t,x)
\]
is a solution of 
\begin{equation}
\label{e5.1}
v_t-v_{xx}- \left( \frac{3}{2t} - \frac{3 \sqrt{\pi}}{2 t^{3/2}} \right)(v-v_x)+e^{-x}v^2=0, \quad \quad x \in \Rm, \quad t > 1.
\end{equation}
We have shown already that $V_{app}$ defined by (\ref{apr724}) is an approximate solution, and the convergence 
\[
u(t,x) \to \phi_\infty(x - x_\infty)
\]
implies that 
\[
|v(t,x) - V_{app}(t,x)| \to 0.
\]
Theorem \ref{t2.2} is an immediate consequence of the definition of $V_{app}$ and the following bound on the error between $v$ and $V_{app}$:
\begin{thm}
\label{t5.1} Given $\gamma > 0$ small, let $V_{app}(t,x)$  be 
the approximate solution constructed in Section~\ref{sec:approx}.
There is $C_\gamma>0$ such that, 
for all $(t,x)\in [1,\infty)\times\R$, we have
\begin{equation}
\label{e5.1001}
\vert e^{x}u(t,x)-V_{app}(t,x)\vert\leq \frac{C_{\gamma}(1+|x|)}{t^{1- \gamma}}.
\end{equation}
\end{thm}

Corollary~\ref{c2.2} also follows from Theorem~\ref{t5.1}.
Let us fix $s\in(0,1)$, let $\sigma_s(t)$ be defined by 
\[
\sigma_s(t)=\sup\{x:~u(t,x)=s\},
\]
and set $\sigma_*^s = \phi_*^{-1}(s)$, so that $\phi_*(\sigma_*^s)=s$. From (\ref{e5.1001}) and the definition of $V^-$, we then have:
\begin{equation}\label{oct736}
\sigma_*^s=\sigma_s(t)+O(t^{-1+\gamma}),
\end{equation}
which is the claim of Corollary~\ref{c2.2} in this moving frame.

\subsection*{The proof of Theorem~\ref{t5.1}}

This is the most technical part of the paper, although the 
idea is really to apply a simple stability argument. 
We will use the self-similar variables 
\begin{equation}\label{oct1014}
\tau=\ln \,t,~~\eta=\farc{x}{\sqrt{t}}
\end{equation}
most of the time.  
As we have noted, there, one may easily reduce the equation for $v$ to an equation on a half-line 
$\eta>0$, 
due to the very fast decay of $v$ for $\eta<0$. Then, we are left with an equation 
for $\eta>0$
that is almost linear: it is perturbed by 
a nonlinear term whose support in~$\eta$ is essentially of the size $e^{-\tau/2}$. 
Moreover, we already know that $e^{-\tau/2} v(\tau,\eta)$ is equivalent, 
for large $\tau$, to 
\[
\alpha_\infty\eta_+e^{-\eta^2/4}.
\]
However, the nonlinear term may be quite large in the small region $\eta\sim O(e^{-\tau/2})$.
We use a bootstrap argument to show that it is in fact harmless, thus opening the way to a classical Liapounov-Schmidt argument of the type~\cite{Sat}.

\subsubsection*{Reduction to the Dirichlet problem}

In view of (\ref{e4.30}), the difference 
$$
\tilde W(t,x)=v(t,x)-V_{app}(t,x).
$$
satisfies an equation
\begin{equation}
\label{e5.8}
\tilde W_t-\tilde W_{xx}- \left(\frac{3}{2t} - \frac{3 \sqrt{\pi}}{2 t^{3/2}}\right)(\tilde W-\tilde W_x)+
e^{-x}(v+V_{app})\tilde W=\tilde E_1(t,x),
\end{equation}
with a function $\tilde E_1$ satisfying:
\begin{equation}\label{oct1106}
|\tilde E_1(t,x)|\le Ct^{-1+3\gamma}
(\un_{0<x<t^\gamma}+e^x \un_{x < 0})+
Ct^{-3/2}e^{-x^2/((4+\gamma)t)}\un_{x>t^\gamma}
+ Ct^{-1+2\gamma}\delta(x-t^\gamma).
\end{equation}
In order to reduce the equation for $\tilde W$ to a Dirichlet problem in the self-similar variables, 
we proceed in several steps.

We first switch to 
\[
W_1(t,x)=\tilde W(t,x)-\tilde W(t,-t^\gamma)\psi(x+t^\gamma).
\]
Here, $\psi(x)$ is a nonegative $C^\infty$ function so that $\psi(x)=1$ for 
$0\leq x\leq1/2$, and $\psi(x)=0$ for $x\geq 1$, so that now $W_1(t,-t^\gamma) = 0$. 
This generates an additional term in the right side 
of \eqref{e5.8} that we denote by $\tilde E_2(t,x)$. 
Taking into account that 
\begin{equation}
\label{e5.20}
v(t,x)+V_{app}(t,x)=O(e^{x})\hbox { for $x<0$}, 
\end{equation}
we obtain
\begin{equation}
\label{e5.9}
\vert \tilde E_2(t,x)\vert\leq Ce^{-t^\gamma}\un_{[0,1]}(x+t^\gamma).
\end{equation}
Next, we translate the origin to $x=-t^{\gamma}$:
the function
\begin{equation}\label{oct820}
W(t,x)=W_1(t,x-t^\gamma) = \tilde W(t, x - t^\gamma) - \tilde W(t,-t^\gamma)\psi(x)
\end{equation}
satisfies
\begin{equation}\label{oct822}
W_t-W_{xx}+\Big(\frac\gamma{t^{1-\gamma}}+\frac3{2t} - \frac{3 \sqrt{\pi}}{2 t^{3/2}}\Big)W_x-
\left(\frac{3}{2t} - \frac{3 \sqrt{\pi}}{2 t^{3/2}} \right) W+e^{t^\gamma-x}(\tilde v+\tilde V_{app})W=G_1(t,x)+
G_2(t,x)
\end{equation}
for $x>0$, with the Dirichlet condition $W(t,0)=0$.
Here, we have introduced
\begin{equation}\label{oct824}
\tilde v(t,x)=v(t,x-t^\gamma),~~\tilde V_{app}(t,x)=V_{app}(t,x-t^\gamma).
\end{equation}
The functions $G_1(t.x)$ and $G_2(t,x)$ in (\ref{oct822}) satisfy
\begin{eqnarray}\label{oct1010}
&&|G_1(t,x)|=|\tilde E_1(t,x-t^\gamma)|  \le  Ct^{-1+3\gamma}
(\un_{t^\gamma<x<2t^\gamma}(x)+e^{x-t^\gamma} \un_{x<t^\gamma}(x))\\
&& ~~~~~~~~~~~~
+ Ct^{-3/2}e^{-(x-t^\gamma)^2/((4+\gamma)t)}\un_{x>2t^\gamma} +Ct^{-1+2\gamma}\delta(x-2t^\gamma),\nonumber
\end{eqnarray}
and
\begin{equation}
\label{oct1012}
\vert G_2(t,x) \vert = | \tilde E_2(t,x - t^\gamma)| \leq Ce^{-t^\gamma}\un_{[0,1]}(x).
\end{equation}

We now express (\ref{oct822}) in the self-similar variables \eqref{oct1014}. With $L$ defined by (\ref{e3.1001}), this gives
\begin{eqnarray}\label{oct826}
W_\tau+\Big(L-\farc{1}{2}\Big)W + e^{\tau+e^{\gamma\tau}-\eta e^{\tau/2}}(\tilde v+\tilde V_{app})W & = & -\Big(\gamma e^{\gamma\tau}+
\farc{3}{2} - \frac{3 \sqrt{\pi}}{2}e^{-\tau} \Big)e^{-\tau/2}W_\eta \no \\
&& - \frac{3 \sqrt{\pi}}{2}e^{-\tau/2} W  + e^{-\eta^2/8} (E_1+E_2),
\end{eqnarray}
with $E_1(\tau,\eta)$ satisfying 
\begin{eqnarray}\label{oct1112}
&&|E_1(\tau,\eta)|\le Ce^{2\gamma\tau}e^{\eta^2/8}
\un\Big(e^{-(\frac{1}{2} - \gamma) \tau} < \eta<2e^{-(\frac{1}{2} - \gamma) \tau}\Big)
\nonumber\\
&&~~~~~~~~~~~~
+Ce^{2\gamma\tau}e^{\eta^2/8} e^{\eta e^{\tau/2} - e^{\gamma\tau}} \un\Big(0<\eta<e^{-(\frac{1}{2} - \gamma) \tau}\Big)
\\
&&~~~~~~~~~~~~
+Ce^{-\tau/2}e^{\eta^2/8}
e^{-(\eta -e^{(-1/2+\gamma)\tau})^2/(4+\gamma)}
\un\Big(\eta>2e^{(-1/2+\gamma)\tau}\Big)\nonumber\\
&&~~~~~~~~~~~~
+Ce^{-(1/2 - 2\gamma) \tau}e^{\eta^2/8}
\delta(\eta-2e^{(-1/2+\gamma)\tau}=
E_{11}+E_{12}+E_{13}+E_{14},\nonumber
\end{eqnarray}
and 
\begin{equation}\label{oct1002}
|E_2(\tau,\eta)|\le e^{\eta^2/8} 
e^{\tau}e^{-e^{\tau\gamma}}\un\Big(0 <\eta < e^{-\tau/2}\Big).
\end{equation}
Notice that the support of $E_{11}$, $E_{12}$, $E_{14}$ is very small, despite the larger prefactor, compared to~$E_{13}$ and $E_2$. Also notice that, in the expression of the Dirac masses, we gain a factor $e^{-\tau/2}$, due to the relation
$$
\delta(x-2\tau^\gamma)=e^{-\tau/2}\delta(\eta-2e^{(-1/2+\gamma)\tau}).
$$
Finally, we symmetrize the operator $L$ by introducing the function
\begin{equation}
\label{e5.7}
w(\tau,\eta)=e^{\eta^2/8}W(\tau,\eta),
\end{equation}
which satisfies
\begin{eqnarray}
\label{e5.10}
&&w_\tau+\M w+
e^{\tau+(\eta_\gamma(\tau)-\eta)e^{\tau/2}}(\tilde v+\tilde V_{app})w= 
\di\sum_{i=1}^2 E_i(\tau,\eta)+E_3(\tau,\eta), \quad \quad \eta > 0
\end{eqnarray}
with the Dirichlet boundary condition $w(\tau,0)=0$.
Here we have defined the operator
\begin{equation}
\label{e5.6}
\M w=- w_{\eta \eta} +\biggl(\frac{\eta^2}{16}-\frac{5}{4}\biggl)w,
\end{equation}
and set
\begin{equation}
\label{e5.14}
\eta_\gamma(\tau)=e^{-(\frac12-\gamma)\tau},~~~~
E_3(\tau,\eta)=-\biggl(\gamma e^{-(\frac12-\gamma)\tau}+\frac{3e^{-\tau/2}}{2}  - \frac{3 \sqrt{\pi}}{2}e^{-3\tau/2} \biggl)(w_\eta-\frac\eta4w) - \frac{3\sqrt{\pi}}{2} e^{-\tau/2} w.
\end{equation}
Strictly speaking, $E_3$ depends on $w$ and $w_\eta$, 
but we omit this dependence for the notational purposes.

Recall that, in the self-similar variables, $V_{app}$ grows as $e^{\tau/2}$. From the convergence result of \cite{NRR1} (Lemma 5.1, in particular) and the definition of $V_{app}$ it follows that
\begin{equation}
\label{e5.5}
\lim_{\tau\to+\infty}e^{-\tau/2}\Vert w(\tau,.)\Vert_{L^2(\R_+)}=0.
\end{equation}
Our goal is to improve this $o(e^{\tau/2})$ bound on $w$ to an exponentially decaying estimate for $w$.

\subsubsection*{From $o(e^{\tau/2})$ to $O(e^{10\gamma\tau})$ asymptotics for the $L^2$ norm of $w$}

The principal eigenfunction of the self-adjoint operator $\M$ with the Dirichlet
boundary condition at~$\eta=0$ is
\begin{equation}
\label{e5.11}
e_0(\eta)=c_0\eta e^{-\eta^2/8}, \quad \quad \quad \M e_0=-\di\frac{e_0}2,
\end{equation}
with the constant $c_0$ chosen so that $\Vert e_0\Vert_{L^2(\R_+)}=1$.   The next eigenvalue is $\lambda_1 = 1/2$ with eignfunction $e_1(\eta) = c_1 e^{\eta^2/8} ( \eta e^{-\eta^2/4})''$; higher eigenfunctions of $\mathcal{M}$ can be expressed in terms of Hermite polynomials. We decompose the solution of~(\ref{e5.10}) as
\begin{equation}
\label{e5.12}
w(\tau)=\langle e_0,w(\tau)\rangle e_0+w^\bot(\tau), \quad \quad \quad 
\int_{\R_+}e_0(\eta)w^\bot(\tau,\eta)d\eta=0.
\end{equation}

{\bf Step 1: a bound for $\langle e_0,w\rangle $.}
We have, projecting \eqref{e5.10} onto $e_0$ and using~\eqref{e5.12}:
\begin{equation}
\label{e5.15}
\frac{d\langle e_0,w\rangle }{d\tau}-\frac{\langle e_0,w\rangle }2+
\langle e_0,e^{\tau+(\eta_\gamma(\tau)-\eta)e^{\tau/2}}(\tilde v+\tilde V_{app})w\rangle 
= \di\sum_{i=1}^3\langle e_0,E_i(\tau)\rangle .
\end{equation}
Let us bound the various perturbative terms in (\ref{e5.15}).  The terms involving $E_1$ and $E_2$ in the right side 
are easily treated. In view of \eqref{oct1112} 
we have
\begin{equation}
\label{oct902}
\vert\langle e_0,E_{1}(\tau)\rangle \vert
\leq Ce^{-(\frac12-3\gamma)\tau}.
\end{equation}
and \eqref{oct1002} implies
\begin{equation}
\label{e5.16}
\vert \langle e_0,E_2(\tau)\rangle \vert \leq Ce^{-e^{\gamma \tau}} \leq Ce^{-(\frac12-3\gamma)\tau},
\end{equation}
as well. As for the term involving $E_3$, using (\ref{e5.14}) and
integrating by parts, we get
\begin{equation}
\label{e5.18}
\vert \langle e_0,E_3(\tau)\rangle \vert \leq
\biggl(\gamma e^{-(\frac12-\gamma)\tau}+\frac{9 e^{-\tau/2}}2\biggl)\left( | 
\langle e_0',w\rangle| + |\langle e_0,\frac\eta4 w\rangle| + | \langle e_0, w\rangle | \right).
\end{equation}
Because of \eqref{e5.5}, we obtain 
\begin{equation}
\label{e5.19}
\vert \langle e_0,E_3(\tau)\rangle \vert\leq Ce^{2\gamma\tau}.
\end{equation}

It finally remains to estimate  the last term in the left side 
of (\ref{e5.15}), and some care should be given to it: 
although the exponential term 
is small outside of the very small set  $0<\eta<\eta_\gamma$,
it could be very large (of the order $e^\tau$) there. This will 
be compensated by the smallness of the factor~$v+V_{app}$.
Let us recall \eqref{e5.20} and (\ref{oct824}) which imply that in the self-similar variables
\begin{equation}
\label{oct906}
\vert \tilde v(\tau,\eta)+ \tilde V_{app}(\tau,\eta)\vert,\ \vert w(\tau,\eta)\vert\leq 
C e^{\eta e^{\tau/2}-e^{\gamma\tau}}=Ce^{e^{\tau/2}(\eta-\eta_\gamma(\tau))}
\hbox{ for $0\leq\eta\leq\eta_\gamma(\tau)$.}
\end{equation}
Let us decompose the inner product
\begin{equation}\label{oct1134}
Q(\tau)=\langle e_0,e^{\tau+(\eta_\gamma(\tau)-\eta)e^{\tau/2}}(\tilde v+\tilde V_{app})w\rangle=
\int_0^{\eta_\gamma(\tau)}+\int_{\eta_\gamma(\tau)}^{\infty} = I_1 + I_2 ,
\end{equation}
For $\eta\leq\eta_{\gamma}(\tau)$ we use the bound  $0 \leq e_0(\eta)\leq c_0\eta$. 
Using \eqref{oct906},
we obtain
\begin{eqnarray}\label{oct912}
&&I_1\le\int_0^{\eta_\gamma(\tau)}e_0(\eta)
e^{\tau+(\eta_\gamma(\tau)-\eta)e^{\tau/2}}(\tilde v+\tilde V_{app})\vert w\vert d\eta
\leq C\int_0^{\eta_\gamma(\tau)}\eta e^{\tau+(\eta-\eta_\gamma(\tau))
e^{\tau/2}}d\eta\nonumber\\
&&~~\leq 
C\eta_\gamma(\tau)e^{\tau}e^{-\tau/2}
=Ce^{\gamma\tau}. 
\end{eqnarray}
As for $I_2$,
we have that
$$
\vert \tilde v(\tau,\eta)+\tilde V_{app}(\tau,\eta)\vert,\ \vert w(\tau,\eta)\vert
\leq C(1+\eta e^{\tau/2})
$$
for all $\eta \in \Rm$. This implies
\begin{eqnarray}\label{oct914}
I_2 & \leq & \int_{\eta_{\gamma}(\tau)}^{\infty} e_0(\eta)e^{\tau+(\eta_\gamma(\tau)-\eta)e^{\tau/2}}(\tilde v+\tilde V_{app})
\vert w\vert d\eta \leq C\int_{\eta_\gamma(\tau)}^\infty \eta e^{\tau+(\eta_\gamma(\tau)-\eta)e^{\tau/2}}(1+\eta e^{\tau/2})^2d\eta
\nonumber\\
& \leq & C e^{2 \tau} \int_{\eta_\gamma(\tau) }^\infty \eta^3  e^{- (\eta - \eta_\gamma(\tau) ) e^{\tau/2} }\, d \eta \leq C(\eta_{\gamma}(\tau))^3 e^{3 \tau/2 }\le Ce^{3 \gamma\tau},
\end{eqnarray}
and therefore,
\begin{equation}\label{oct1136}
|Q(\tau)|\le C e^{3\gamma\tau}.
\end{equation}

Putting everything together, we infer 
that
\begin{equation}\label{oct918}
\frac{d\langle e_0,w\rangle }{d\tau}-\frac{\langle e_0,w\rangle }2= \varphi(\tau),
\end{equation}
with
\[
|\varphi(\tau)|\le Ce^{3\gamma\tau}.
\] 
We see that
\begin{equation}\label{oct922}
\farc{d}{dt}\Big(\langle e_0,w\rangle e^{-\tau/2}\Big)=
\varphi(\tau)e^{-\tau/2}.
\end{equation}
Taking into account \eqref{e5.5}, we can 
integrate (\ref{oct922}) 
from $\tau$ to $+\infty$ leading to
\begin{eqnarray}\label{oct920}
\langle e_0,w(\tau)\rangle=
-\int_\tau^{+\infty}e^{(\tau-\tau')/2}\varphi(\tau')d\tau',
\end{eqnarray}
hence
\begin{equation}\label{oct926}
|\langle e_0,w(\tau)\rangle|\le 
C\int_\tau^{+\infty}e^{(\tau-\tau')/2}e^{3\gamma\tau '}d\tau'\le 
C_\gamma e^{3\gamma\tau}.
\end{equation}
This bound will be improved in the next step.

{\bf Step 2. An $L^2$ bound for $w^\bot(\tau)$.} We multiply (\eqref{e5.10}) by $w^\bot$, and integrate by parts: 
\begin{equation}
\label{e5.24}
\frac12\frac{d\Vert w^\bot\Vert^2}{d\tau}+
\langle\M w^\bot,w^\bot\rangle
+\int_{\R_+}e^{\tau+(\eta_\gamma(\tau)-\eta)e^{\tau/2}}(\tilde v
+\tilde V_{app})ww^\bot d\eta=\sum_{i=1}^3\int_{\R_+}E_iw^\bot d\eta.
\end{equation}
We denoted here the $L^2(\R_+)$ norm by $\Vert\cdot\Vert$. 
Once again, we need to bound the perturbative terms in~(\ref{e5.24}). Let us start with the less standard term:
\[
q(w):= \di\int_{\R_+}e^{\tau+(\eta_\gamma(\tau)-\eta)e^{\tau/2}}
(\tilde v+\tilde V_{app})ww^\bot d\eta=J_1(\tau)+J_2(\tau),
\]
with the two terms coming from 
the decomposition  \eqref{e5.12} for $w$. We have
\begin{equation}\label{oct1178}
J_1(\tau)=\langle e_0,w(\tau)\rangle \int_{\R_+}e^{\tau+(\eta_\gamma(\tau)-\eta)e^{\tau/2}}(\tilde v+ \tilde V_{app})e_0w^\bot d\eta.
\end{equation}
We know from Step 1 that 
\[
\langle e_0,e^{\tau+(\eta_\gamma(\tau)-\eta)e^{\tau/2}}(\tilde v+ \tilde V_{app}) |w|\rangle = |Q(\tau)| \le Ce^{3\gamma\tau}.
\]
Together with (\ref{oct926}) this gives
\begin{equation}\label{oct930}
|J_1(\tau)|\le Ce^{6\gamma\tau}.
\end{equation}
Furthermore, $J_2(\tau)$ is positive, so we do not need to estimate it.

As for the three terms in the right side of (\ref{e5.24}),
in view of \eqref{oct1112}  we have, with some
constant~$C_\gamma>0$: first,
\begin{equation}\label{oct1020}
|\langle w^\bot,E_{11}\rangle|+
|\langle w^\bot,E_{12}\rangle|\le \gamma \|w^{\bot}\|^2+
\farc{1}{4 \gamma}(\|E_{11}\|^2+\|E_{12}\|^2)\le \gamma \|w^{\bot}\|^2+
C_\gamma e^{5\gamma\tau}e^{-\tau/2},
\end{equation}
while for $E_{13}$ we have
\begin{equation}\label{oct1028}
|\langle w^\bot,E_{13}\rangle|\le \gamma \|w^{\bot}\|^2+
\farc{1}{4\gamma}\|E_{13}\|^2\le \gamma \|w^{\bot}\|^2+
\farc{C}{\gamma}e^{-\tau}.
\end{equation}
Finally, for $E_{14}$ we have
\begin{eqnarray}\label{oct1020bis}
|\langle w^\bot,E_{14}\rangle| &\le & Ce^{2\gamma\tau}
|w^\bot(\tau,2e^{-\tau/2+\gamma\tau})|\le 
Ce^{2 \gamma \tau} e^{-\frac{\tau}{4} + \frac{\gamma}{2}\tau}\| \partial_\eta w^\bot(\tau,\cdot)\|_{L^2}\\
&\leq & Ce^{2 \gamma \tau} e^{-\frac{\tau}{4} + \frac{\gamma}{2}\tau}(1+\langle Mw^\bot,w^\bot\rangle+\|w^\bot\|^2).
\nonumber
\end{eqnarray}
For $E_2$ we may simply estimate
\begin{equation}\label{oct1122}
|\langle w^\bot,E_{2}\rangle| \le    \gamma\|w^\bot\|^2 + \frac{1}{4\gamma} \norm{E_2}^2 \leq \gamma\|w^\bot\|^2+
C_\gamma e^{2\tau-2e^{\gamma\tau}}e^{-\tau/2}\le
\gamma\|w^\bot\|^2+
C_\gamma e^{-\tau/2}.
\end{equation}
As for $E_3$, we have
\begin{equation}\label{oct1114}
\Big\vert\int_{\R_+}(w_\eta-\frac\eta4w)w^\bot d\eta\Big|
\leq \int \eta^2 (w^\bot)^2d\eta+
C_\gamma\langle e_0,w\rangle ^2+\gamma\Vert w^\bot\Vert^2\le
C\|w^\bot\|^2+C\langle\M w^\bot,w^\bot\rangle+
C_\gamma\langle e_0,w\rangle ^2,
\end{equation}
hence
\begin{equation}\label{oct1116}
|\langle w^\bot,E_{3}\rangle|\le C_\gamma e^{(-1/2+\gamma)\tau}
\Big(\|w^\bot\|^2+\langle\M w^\bot,w^\bot\rangle+
\langle e_0,w\rangle ^2\Big).
\end{equation}
Recall that the second eigenvalue of $\mathcal{M}$ is $1/2$, so we have
\[
\langle \M w^\bot, w^\bot\rangle
\geq\di\frac{\Vert w^\bot \Vert^2}2.
\]
Putting everything together, this yields
\begin{equation}\label{oct1132}
\frac12\frac{d\Vert w^\bot\Vert^2}{d\tau}+
\biggl(\frac12-\gamma-C_\gamma e^{-(\frac{1}{4}-\frac{3\gamma}{2})\tau}\biggl)
\Vert w^\bot\Vert^2\leq |J_1(\tau)|\le Ce^{6\gamma\tau}.
\end{equation}
This implies 
\begin{equation}
\label{e5.25}
\Vert w^\bot\Vert\leq C_\gamma e^{3\gamma\tau}.
\end{equation}
Because of (\ref{oct926}), this bound also holds for the full solution: $\norm{w} \leq C_\gamma e^{3\gamma\tau}$. 

\subsubsection*{Upgrading the $L^2$ bound for $w$ to an $L^\infty$  bound.} We now know that $w$ satisfies a linear inhomogeneous equation of the form
\begin{align}
w_\tau + \mathcal{M} w + H(\tau,\eta) w + g(\tau)(w_\eta - \frac{\eta}{4} w) + \frac{3 \sqrt{\pi}}{2} e^{-\tau/2} w & = f(\tau,\eta) + O(e^{2 \gamma\tau}) \un_{[0,2 \eta_\gamma(\tau))} \nonumber \\
& \quad + h(\tau) \delta(\eta - 2 \eta_\gamma(\tau)) \label{weqnnew}
\end{align}
with $w(\tau,0) = 0$, where 
\begin{equation}
\label{e5.85}
g(\tau) = \biggl(\gamma e^{-(\frac12-\gamma)\tau}
+\frac{3e^{-\tau/2}}2 - \frac{3 \sqrt{\pi}}{2}e^{-3 \tau/2} \biggl),
\end{equation}
and 
\[
H(\tau,\eta) = e^{\tau+(\eta_\gamma(\tau)-\eta)e^{\tau/2}}(\tilde v+\tilde V_{app}) \geq 0.
\]
The forcing terms $f$ and $h$ satisfy
\[
|f(\tau,\eta)| \leq C e^{-\tau/2} e^{- \eta^2/16},
\]
and 
\begin{equation}
\label{e5.86}
h(\tau), h'(\tau) = O(e^{-(1/2-2 \gamma) \tau}).
\end{equation}
For the moment we are not going to use the full force of this estimate, we will only use the fact that $h$ and $h'$ grow at most like $e^{2\gamma\tau}$.
Notice that, for every $a>0$, the singular term on the right side of \eqref{weqnnew} is supported in $[0,a/2)$ for $\tau$ large enough. Also, for every $a > 0$,
$$
\lim_{\tau\to+\infty}\Vert H(\tau,.)\Vert_{L^\infty((a,+\infty)}=0.
$$
Hence, by parabolic regularity (e.g. \cite{Lieb}, Theorem 6.30, 7.43) and the bound $\norm{w}_{L^2} \leq C_\gamma e^{3\gamma\tau}$,  we infer that
\[
\Vert w \Vert_{L^\infty([a,A])}\leq C_{a,A}e^{5\gamma\tau},
\]
for $a$ small, $A$ large. The $L^\infty$ estimates on the perturbative terms in 
the equation (\ref{e5.10}) for $w$ imply that 
for $\eta\geq A$ sufficiently large, $w(\tau,\eta)$ cannot attain its maximum
at a point $\eta>A$ where it is larger than $Ce^{5\gamma\tau}$, thus we
have
\begin{equation}
\label{e5.26}
\Vert w\Vert_{L^2(\R_+)}+\Vert w\Vert_{L^\infty((a,+\infty))}\leq C_{a,\gamma} e^{10\gamma\tau},
\end{equation}
for $a>0$ small.
 To retrieve the $L^\infty$ bound on the full half line, we proceed as follows. 
By the Kato inequality, equation~\eqref{weqnnew} for~$w$   
yields, writing out explicitly the operator $\M$:
\begin{eqnarray}
\label{e5.27}
&&\!\!\!\!\!\!\!\!\!\!\!\!\!\!
\partial_\tau\vert w\vert-|w|_{\eta\eta}+\Big(\farc{\eta^2}{16}-\frac{5}{4}\Big)
|w|
+ g(\tau) \biggl(\partial_\eta\vert w\vert-\frac{\eta}{4} \vert w\vert\biggl) \no \\
&& \quad \leq Ce^{-(\frac12-\gamma)\tau} +Ce^{2\gamma\tau}
\un\Big(0<\eta<2 \eta_\gamma(\tau) \Big)+
Ce^{2\gamma \tau}
\delta(\eta-2 \eta_\gamma(\tau)),
\end{eqnarray}
with $g(\tau)$ given by \eqref{e5.85}.
Let $a\in(0,1)$ be small enough so that (\ref{e5.27}) implies
\begin{eqnarray}
\label{oct1140}
\partial_\tau\vert w\vert-|w|_{\eta\eta} -10|w|
+g(\tau)\partial_\eta\vert w\vert \leq Ce^{2\gamma\tau}+Ce^{2\gamma \tau}
\delta(\eta-2 \eta_\gamma(\tau)),
\end{eqnarray}
for $\eta \in (0,a)$ with the boundary conditions 
\begin{equation}\label{oct1142}
|w|(\tau,0)=0,~~|w|(\tau,a)\le C_{a,\gamma}e^{10\gamma\tau},
\end{equation}
which is achievable, due to (\ref{e5.26}). Drop the subscript ${}_{a,\gamma}$ - it is not useful anymore here - and let us write
\[
|w|(\tau,\eta)\le Ce^{10\gamma\tau}\psi(\tau,\eta)+e^{2\gamma \tau}
\phi(\tau,\eta),
\]
with the function $\psi(\tau,\eta) \geq 0$ such that
\begin{eqnarray}
\label{oct1144}
&&\partial_\tau\psi-\psi_{\eta\eta} -11\psi
+g(\tau)\partial_\eta\psi
=Ce^{-8\gamma\tau},\\
&&\psi(\tau,0)=0,~~\psi(\tau,a)=1.\nonumber
\end{eqnarray}
Possibly decreasing $a$, we may ensure that
the principal eigenvalue $\lambda_a$ of the 
Dirichlet Laplacian on the interval~$(0,2a)$ is 
sufficiently large, say, $\lambda_a>100$. Then there exists a constant
$C>0$ so that 
\begin{equation}\label{oct1146}
\psi(\tau,\eta)\le C\eta.
\end{equation}
We choose the function $\phi \geq 0$ so that it satisfies
\begin{eqnarray}
\label{oct1170}
\partial_\tau\phi-\phi_{\eta\eta} -11\phi
+g(\tau)\partial_\eta\phi =C \delta(\eta-2e^{(-1/2+\gamma)\tau}),
\end{eqnarray}
with the boundary conditions 
\begin{equation}\label{oct1172}
\phi(\tau,0)=0,~~\phi(\tau,a)=0.
\end{equation}
Let us prove that
\begin{equation}\label{oct1174}
\phi(\tau,\eta)\le C\eta.
\end{equation}
We have $\phi(\tau,\eta)=\phi_0(\tau,\eta)+\phi_1(\tau,\eta)$
with
\begin{eqnarray*}
&&-\partial_{\eta\eta}\phi_0=C 
\delta(\eta-2e^{(-1/2+\gamma)\tau})\\
&&\phi_0(\tau,0)=\phi_0(\tau,a)=0,
\end{eqnarray*}
and
\[
\partial_\tau\phi_1-\partial_{\eta\eta}\phi_1 -11\phi_1
+g(\tau)\partial_\eta\phi_1
=-\partial_\tau\phi_0-g(\tau)\partial_\eta\phi_0+11\phi_0,
\]
with the boundary conditions
\[
\phi_1(\tau,0)=0,~~\phi_1(\tau,a)=0.
\]
The function $\phi_0$ is easily computed:
\br
\phi_0(\tau,\eta)= \left\{\begin{array}{rll} \di C \frac{(a-\xi_\gamma(\tau))}{a}\eta,\   \  &\eta \leq\xi_\gamma(\tau) \\ 
\di C \frac{(a-\eta)\xi_\gamma(\tau)}{a},\   \  &\eta \geq\xi_\gamma(\tau)
\end{array}
\right.
\er
with $\xi_\gamma(\tau)=2e^{(-1/2+\gamma)\tau}$.  So, all the quantities 
$\phi_0$, $\partial_\eta\phi_0$ and $\partial_\tau\phi_0$ are uniformly bounded, hence 
(recall $\lambda_a\geq 100$) we have \eqref{oct1174}.
It follows that
\begin{equation}
\label{e5.29}
\vert w(\tau,\eta)\vert\leq Ce^{10\gamma\tau}\eta\   
\hbox{ for $\tau\geq 0$ and $0\leq\eta\leq a$.}
\end{equation}
This not only yields the full $L^\infty$ estimate for $w$, this gives an extra information on how $w(\tau,\eta)$ grows in the vicinity of 0,
that we are going to use in our next step.

\subsubsection*{From the $O(e^{10\gamma\tau})$ growth to 
$O(e^{-(\frac12-100\gamma)\tau})$ decay for $\Vert w\Vert_{L^2}$} 

The next step is thus to improve the ``slow" $O(e^{10\gamma\tau})$ growth in (\ref{e5.26})
to actual decay in time. Let us come back to (\ref{e5.15}), the equation
for $\langle e_0,w\rangle$:
\begin{equation}
\label{oct1126}
\frac{d\langle e_0,w\rangle }{d\tau}-\frac{\langle e_0,w\rangle }2+
\langle e_0,e^{\tau+(\eta_\gamma(\tau)-\eta)e^{\tau/2}}(\tilde v+\tilde V_{app})w\rangle 
= \di\sum_{i=1}^3\langle e_0,E_i(\tau)\rangle .
\end{equation}
The bounds (\ref{oct902}) and (\ref{e5.16}) are already of the 
``good" size $O(e^{-(1/2-3\gamma)\tau})$, and the already obtained bound
(\ref{e5.26}) allows us to improve (\ref{e5.19}) to
\begin{equation}
\label{oct1130}
\vert \langle e_0,E_3(\tau)\rangle \vert \leq 
\biggl(\gamma e^{-(\frac12-\gamma)\tau}+\frac{9e^{-\tau/2}}2\biggl) \left( | 
\langle e_0',w\rangle| + |\langle e_0,\frac{\eta}{4} w\rangle| + |\langle e_0, w \rangle|\right) \leq
C e^{(-1/2+15\gamma)\tau}.
\end{equation}
Thus, what really limits the decay improvement for $\langle e_0,w\rangle$
is the integral 
\begin{equation}\label{oct1150}
Q(\tau)=\langle e_0,e^{\tau+(\eta_\gamma(\tau)-\eta)e^{\tau/2}}(\tilde v+\tilde V_{app})w\rangle ,
\end{equation}
that we have so far only managed to bound by $Ce^{3\gamma\tau}$ (see 
(\ref{oct1136})). We have already noted that the integrand could 
be very large only for $\eta$ of the order 
\[
\eta_\gamma(\tau)=e^{(-1/2+\gamma)\tau}.
\]
On the other hand, from \eqref{e5.29}, $w$ has a bounded linear growth
in a neighborhood of $\eta=0$. This will bring a 
small factor of the order $\eta$ in the integrand, which will, in 
turn, make the integral be of a smaller order. 
 
So, let us consider $Q(\tau)$ given by (\ref{oct1150}). Using (\ref{oct906}) and
(\ref{e5.29}), we deduce the following improvement of~(\ref{oct912}):
\begin{eqnarray}\label{oct1160}
&&I_1\le\int_0^{\eta_\gamma(\tau)}e_0(\eta)
e^{\tau+(\eta_\gamma(\tau)-\eta)e^{\tau/2}}(\tilde v+ \tilde V_{app})\vert w\vert d\eta
\leq Ce^{10\gamma\tau}\int_0^{\eta_\gamma(\tau)}\eta^2 
e^{\tau}d\eta\nonumber
\\
&&~~~\leq 
Ce^{10\gamma\tau}[\eta_\gamma(\tau)]^3e^{\tau} =Ce^{(-1/2+20\gamma)\tau},
\end{eqnarray}
while (\ref{oct914}) can be improved to
\begin{eqnarray}\label{oct1162}
I_2 & \leq &  \int^{a}_{\eta_\gamma(\tau)}
e_0(\eta)e^{\tau+(\eta_\gamma(\tau)-\eta)e^{\tau/2}}(\tilde v+\tilde V_{app})
\vert w\vert d\eta + C e^{3 \tau/2} e^{- a/2 e^{\tau/2}} \no \\
& \leq & C e^{10\gamma\tau}\int^{a}_{\eta_\gamma(\tau)}
\eta e^{\tau+(\eta_\gamma(\tau)-\eta)e^{\tau/2}}(1+\eta e^{\tau/2})\eta d\eta
\nonumber\\
& \leq & C e^{10\gamma\tau} e^{\tau} e^{\tau/2} \int^{a}_{\eta_\gamma(\tau)} 
\eta^3 e^{-(\eta - \eta_\gamma(\tau))e^{\tau/2}}d\eta  \leq Ce^{10\gamma\tau}(\eta_{\gamma}(\tau))^3 e^{\tau}
\le Ce^{(-1/2+20\gamma)\tau}.
\end{eqnarray}
Equation \eqref{oct1126} for $\langle e_0,w(\tau)\rangle $ now gives
\begin{equation}
\label{e5.30}
|\langle e_0,w(\tau)\rangle|\le C\int_\tau^{+\infty}e^{(\tau-\tau')/2}
e^{(-\frac12+20\gamma)\tau'}d\tau'\le Ce^{-(\frac12-20\gamma)\tau}.
\end{equation}

Moreover, equation (\ref{oct1132}) for $w^\bot$
shows that the only ``slightly large" term that 
potentially can make $w^\bot(\tau,\eta)$ grow in $\tau$ 
is $J_1(\tau)$ given by (\ref{oct1178})
\begin{equation}\label{oct1180}
J_1(\tau)=\langle e_0,w(\tau)\rangle \int_{\R_+}e^{\tau+(\eta_\gamma(\tau)-\eta)e^{\tau/2}+x_0}(\tilde v+ \tilde V_{app})e_0w^\bot d\eta.
\end{equation}
However, we may now use (\ref{e5.30}) to bootstrap (\ref{oct930}) to
\begin{equation}\label{oct1182}
|J_1(\tau)|\le Ce^{(-1/2+40\gamma)\tau}.
\end{equation}
Using this in (\ref{oct1132}) gives us
\begin{equation}\label{oct1184}
\Vert w^\bot\Vert\leq Ce^{-(\frac12-50\gamma)\tau}.
\end{equation}
This implies the same estimate for the full solution $w$. As in the passage
from (\ref{e5.25}) to (\ref{e5.26}) we obtain
\begin{equation}
\label{oct1186}
\Vert w\Vert_{L^2(\R_+)}+\Vert w\Vert_\infty\leq C_\gamma e^{(-1/2+100\gamma)\tau}.
\end{equation}

\subsubsection*{Concluding the proof of Theorem \ref{t5.1}} 
The last step seems to yield a $t^{\gamma-1/2}$ decay for $w$. 
However, recall that we want a $t^{\gamma-1}$ estimate. To this end, it suffices to remember that $w(\tau,\eta)$ solves a Dirichlet problem, hence $w$ 
should have an extra $\eta$ factor. To show that, it suffices to argue just as in the proof of estimate \eqref{e5.29}, up to the fact that, this time, the slow $e^{10\gamma\tau}$ growth is replaced by the decay $e^{-(1/2-100\gamma)\tau}$, and that we may use the full estimate \eqref{oct1186}. Repeating this argument, we end up with
\begin{equation}
\label{e4.150}
|w(\tau,\eta)|\leq C_\gamma \eta e^{-(1/2-100\gamma)\tau}.
\end{equation}
To obtain the conclusion of Theorem \ref{t5.1}, it suffices to unzip \eqref{e4.150}, reverting to the $(t,x)$ variables. We obtain
\begin{equation}\label{oct1188}
\vert v(t,x)-V_{app}(t,x)\vert\leq \frac{C}{t^{\frac12-100\gamma}}
\frac{x+t^\gamma}{\sqrt t}, \quad \quad \text{for}\;\;x > - t^\gamma + 2,\;\; t \geq 1.
\end{equation}
This implies Theorem \ref{t5.1}.~$\Box$

\end{document}